\documentclass[12pt]{article}

\usepackage{pdfsync, hyperref}
\usepackage{amssymb, amsfonts}
\usepackage{latexsym}
\usepackage{amsmath}
\usepackage{amscd}
\usepackage{amsthm}
\usepackage{xy}
\xyoption{all}
\CompileMatrices

\usepackage[mathscr]{eucal}

\usepackage[T3,T1]{fontenc}
 \DeclareSymbolFont{tipa}{T3}{cmr}{m}{n}
 \DeclareMathAccent{\invbreve}{\mathalpha}{tipa}{16}
 
  \DeclareFontFamily{OT1}{pzc}{}
  \DeclareFontShape{OT1}{pzc}{m}{it}{<-> s * [1.050] pzcmi7t}{}
\def\mathpzc#1{\mbox{\fontfamily{pzc}\itshape#1}}
\newenvironment{frcseries}{\fontfamily{frc}\selectfont}{}
 \newcommand{\textfrc}[1]{{\frcseries#1}}
\newcommand{\mathfrc}[1]{\text{\textfrc{#1}}}

\DeclareMathAlphabet\eusm{U}{eus}{m}{n}
\def\makebb#1{\expandafter\def\csname bb#1\endcsname{{\mathbb{#1}}}\ignorespaces}
\def\makerm#1{\expandafter\def\csname rm#1\endcsname{{\rm #1}}\ignorespaces}
\def\makebf#1{\expandafter\def\csname bf#1\endcsname{{\bf #1}}\ignorespaces}
\def\makegr#1{\expandafter\def\csname gr#1\endcsname{{\mathfrak{#1}}}\ignorespaces}
\def\makescr#1{\expandafter\def\csname scr#1\endcsname{{\mathscr{#1}}}\ignorespaces}
\def\makecal#1{\expandafter\def\csname cal#1\endcsname{{\cal #1}}\ignorespaces}
\def\makeudl#1{\expandafter\def\csname udl#1\endcsname{{\underline{#1}}}\ignorespaces}
\def\doLetters#1{%
  #1A #1B #1C #1D #1E #1F #1G #1H #1I #1J #1K #1L #1M
  #1N #1O #1P #1Q #1R #1S #1T #1U #1V #1W #1X #1Y #1Z}
\def\doletters#1{%
  #1a #1b #1c #1d #1e #1f #1g #1h #1i #1j #1k #1l #1m
  #1n #1o #1p #1q #1r #1s #1t #1u #1v #1w #1x #1y #1z}
\doLetters\makebb\doLetters\makecal\doLetters\makerm\doLetters\makebf\doLetters\makescr
\doletters\makerm \doletters\makebf \doLetters\makegr
\doletters\makegr \doLetters\makeudl \doletters\makeudl

\usepackage{graphicx}
\usepackage{calc}

\newlength{\wcwidth}
\newlength{\wcheight}
\newcommand{\widecheck}[1]{\ensuremath{
\settowidth{\wcwidth}{#1}
\settoheight{\wcheight}{#1}
\addtolength{\wcheight}{1pt}
\makebox[0cm][l]{%
\raisebox{\depth+\wcheight}[0cm][0cm]{%
\scalebox{-1}{$\widehat{\hphantom{#1}}$}}}#1
\rule{0pt}{\wcheight+2.5pt}}}

\theoremstyle{plain}

\newtheorem{thm}{Theorem}[section]\newtheorem{prop}[thm]{Proposition}
\newtheorem{lemma}[thm]{Lemma}

\newtheorem*{thm*}{Theorem}
\newtheorem*{conj*}{Conjecture}

\theoremstyle{definition}
\newtheorem{defn}{Definition}[section]

\theoremstyle{remark}

\def\BTitem#1\ETitem{\begin{equation}\hss\left\{\mkern-20mu\parbox{0.9\hsize}{%
\begin{itemize}#1\end{itemize}}\right.\end{equation}}

\numberwithin{equation}{section}

 \newcommand{\epf}{\hfill$\Box$\medbreak}

\def\co{\colon\thinspace}
\def\pf{\noindent{\it Proof.}\enspace}

\newcommand{\slp}{\partial\hspace{-1ex}/}
\newcommand{\ts}{\tilde{\sigma }}

\def\ff{\mathpzc{f}}

\newcommand{\op}{\operatorname}

\newcommand{\sym}{\op{Sym}}

\newcommand{\sO}{\textsc{o}}

\newcommand{\hf}{\textsl{HF}}

\def\HM{\mathop{\textit{HM} } \nolimits}
\def\Hf{\mathop{\textit{HF} }\nolimits}

\def\frcs{\mathfrc{s}}
\def\ds{\displaystyle}

 \def\ul{\underline}
\def\td{\tilde}

 \def\Spin{\mathop{\mathrm{Spin}}\nolimits}

\newcommand{\Conn}{\op{Conn}}

\begin{document}

\title{Morse theory and Seiberg-Witten moduli spaces of
3-dimensional cobordisms, I}

 \author{Yi-Jen Lee} 
 
\date{Preliminary}
\maketitle
\begin{abstract} 


Motivated by a variant of Atiyah-Floer conjecture proposed
in \cite{L2} and its potential generalizations, in this article and
its sequel we study as a first
step properties of  moduli spaces of Seiberg-Witten equations on a
3-dimensional cobordism with cylindrical ends (CCE) \(Y\), perturbed by closed 2-forms of the form \(r*d\ff+w\), where
\(r\geq 1\), where \(\ff\) is a harmonic Morse function with certain
linear growth at the ends of \(Y\), and \(w\) is a certain closed
2-form. 
\end{abstract}

\section{Introduction}\label{sec:intro}

\begin{defn}\label{def:CCE}
 A {\em 3-dimensional cobordism with cylindrical ends} (``CCE'' for
 short) is a connected complete oriented Riemannian 3-manifold \(Y\), such that:
 there is a compact  3-dimensional submanifold with boundary, \(Y_c\subset Y\), and an isometry
\[\iota\co  (-\infty, 0)_t\times \Sigma_-\sqcup(0, \infty)_t\times
\Sigma_+\to Y\backslash Y_c,\] where:
\begin{itemize}
\item \(\Sigma _\pm\) are nonempty oriented compact surfaces;
  \item \((-\infty, 0)_t\times \Sigma_-\)
and \((0, \infty)_t\times\Sigma_+\) are both equipped with a  
product metric, with the first factor endowed with the metric induced
from the affine metric on \(\bbR_t\). 
\end{itemize}
A metric satisfying the above constraints is called a {\em cylindrical
  metric}. \(Y\) is called a {\em CCE from \(\Sigma _-\) to \(\Sigma
  _+\)}, and is frequently denoted by \(Y\co \Sigma _-\to \Sigma_+\). 
We call  \(\scrE_-[Y]:=\iota \big((-\infty, -2)_t\times
\Sigma_-\big)\) and \(\scrE_+[Y]:=\iota\big( (2, \infty)_t\times
\Sigma_+)\) respectively the
{\em negative end} and the {\em positive end} of \(Y\). 
\end{defn}
\begin{defn}\label{def:adm-f}
  Let \(Y\) be a CCE from \(\Sigma _-\) to \(\Sigma _+\), and adopt
  the notations from Definition \ref{def:CCE}. An {\em
    admissible function} \(\ff\co Y\to \bbR\) is a harmonic Morse function
 such that:
  \begin{itemize}
  \item[(1)] It has finitely many critical points. Thus, we may and
    will choose to define \(Y_c\) such that all critical points of
    \(\ff\) lie in the interior of \(Y_c\);
   \item[(2)]  There exists constants \(C_\pm\in\bbR\) such that \(\ff-(\iota
     ^*t+C_\pm)\in L^2_1(\scrE_\pm[Y])\) respectively. Here,  where
     \(t \co (-\infty, 0)\times \Sigma _-\sqcup(0,
        \infty)\times \Sigma _+\to \bbR\backslash \{0\}\) denotes the
        projection to the first factor. 
   \end{itemize}
   \end{defn} 

Let \((Y, \grs)\) be a \(\Spin^c\) 3-manifold, and let \(\bbS\) denote
the associated spinor bundle. Let  \(\op{Conn}(\bbS)\) denote the space
of \(\Spin\)-connections on \(\bbS\). Let \(\rho : \bigwedge^*T^*M\to
\op{End}(\bbS)\) denote the Clifford action, with the
convention\footnote{This convention agrees with that in \cite{KM}, but
  is opposite to that in \cite{L1}.}
chosen such that
\[
  \rho (*\nu )=-\rho (\nu ) \qquad \forall
\nu \in \Omega^1(M).\]
The 3-dimensional {\em Seiberg-Witten equation} on
\((Y, \grs)\) concerns an element \((A, \Psi )\in
\op{Conn}(\bbS)\times\Gamma (\bbS)\), called a (Seiberg-Witten) {\em configuration}, and takes the following general form: 
 \begin{equation}\label{eq:SW3}
  \mathfrak{F}_\mu (A, \Psi ):=\left( \begin{array}{c}\frac{1}{2}*F_{A^t}+\rho^{-1}(\Psi\Psi^*)_0+\frac{i}{4}*\mu\\
     \slp_A\Psi 
 \end{array}
\right)=0, 
 \end{equation}
where \(\mu \) is a closed 2-form (the previously mentioned
perturbation form on whose cohomology class the monopole Floer
homology depends on).  \( \rho \) stands for the Clifford
action, and \(\slp_A\) is the Dirac operator. The subscript
 \((\cdot)_0\) in \((\Psi\Psi^*)_0\) means the tracelss part, and
 \(A^t\) is the connection on \(\det \bbS\) induced from \(A\).
 (In
 general, a further abstract perturbation is needed to make the Floer
 homology well-defined, but that is unnecessary in our context.)
Note that \(A^t\in \op{Conn}(\det \bbS)\) together with the
Levi-Civita connection determines a \(\Spin^c\)-connection \(A\); so we
shall use \(A\) and \(A^t\) interchangeably to specify a \(\Spin^c\) connection.

There is an action of \(C^\infty(Y; U(1))\) on 
\(\op{Conn}(\det \bbS)\times\Gamma (\bbS)\) given by
\[
  g\cdot (A^t, \Psi )=(-2g^{-1}dg, g\Psi )\quad g\in C^\infty(Y; U(1)),
\]
called the {\em gauge action}. 
To configurations are said to be {\em gauge equivalent} if they are
related by such an action. Note that \(\grF_\mu \) is invariant under
the gauge action; so we may refer to a gauge equivalence class as a
solution to the Seiberg-Witten equation (\ref{eq:SW3}).

 Let \(Y\) be a CCE from \(\Sigma _-\) to \(\Sigma _+\), and let \(\ff:
 Y\to \bbR\) be an admissible function.  Fix a \(\Spin^c\) structure
 \(\grs\), and write \(c_1(\grs)=c_1(\det\bbS)\), where \(\bbS\) is
 the spinor bundle associated to \(\grs\). We define the {\em degree} of \(\grs\) (relative to
 \(f\)), denoted by \(d=d_\grs\), by the formula
\[c_1(\grs)|_{\Sigma _{\rm min}}=2d+\chi(\Sigma _{\rm min}),\]
 where  \(\Sigma _{\rm min}\subset Y\) is a regular level surface of
 \(\ff\) which has
 minimal genus among all level surfaces of  \(\ff\). Note that the
 preceding definition does not depend on the choice of \(\Sigma _{\rm
   min}\).

 Let \(w\) be an admissible two form on \(Y\), as defined in
 Definition \ref{def:adm-w} below. We consider a family of perturbation forms
 parametrized by \(r\in \bbR\):
 \begin{equation}\label{def:mu-r}
  \mu _r=\mu _{r,w}=r*d\ff+w, \quad r\geq 1.
\end{equation}
The notion of an admissible configurations is introduced in Definition
\ref{def:adm-config} below. We use  \(\scrZ_{r, w}(Y, \grs;\ff)\) to
denote space of gauge equivalence classes of admissible solutions
\((A, \Psi )\) to \(\grF_{\mu _{r,w}}(A, \Psi )=0\).
 
Given a compact K\"ahler surface, we endow the symmetric product \(\sym^k\Sigma \) is equipped
with the K\"ahler form induced by its identification with the moduli
space of vortex solutions on a degree \(k\) line bundle on
\(\Sigma\). (See e.g. \cite{G}.)  
 \begin{thm}\label{thm:main1}
   Let \(Y\co \Sigma _-\to \Sigma _+\) be a CCE, and \(\ff\) is an
   admissible function. Fix a \(\Spin^c\) structure \(\grs\) on \(Y\)
   with degree \(d\), and let \(d_\pm:=d+\frac{\chi(\Sigma _{\rm
       min})-\chi (\Sigma _\pm)}{2}\).  Let \(\hat{\scrW}\) be the
   space of admissible 2-forms defined following Definition
   \ref{def:adm-w}. Then there exists a constant \(r_0\geq 1\)
   depending only on the metric, \(d\) and \(w\), such that 
    \(\forall r\geq r_0\), there is a Baire subset
   \(\hat{\scrW}_{reg}\subset \hat{\scrW}\) such that \(\forall w\in
   \hat{\scrW}_{reg}\), \(\scrZ_{r, w}(Y,
   \grs;\ff)\) is empty when \(d<0\), and otherwise an orientable smooth manifold of dimension \(d_-+d_+\). It is
   equipped with an ``end point map'' \[
     \Pi _{-\infty}\times \Pi _{+\infty}\co \scrZ_{r,
       w}(Y,\grs;\ff)\to (-\scrV_{r, d_-}(\Sigma_-)\times
     \scrV_{r,d_+}(\Sigma _+),\] 
   which is a Lagrangian immersion. In the above, \(\scrV_{r,d}(\Sigma
   )\) denotes the moduli space of the solutions to the version of
   vortex equation defined in (\ref{def:mod-vor}). As explained in
   \cite{G}, \(\scrV_{r,d}(\Sigma)\simeq\sym^{d}\Sigma\) and is
   endowed with a natural symplectic structure. 
 \end{thm}

 
   \subsection{Notation and Conventions}
   \begin{itemize}
     \item Let \(V\to Y\) be a euclidean or hermitian vector bundle
       over a manifold (possibly with boundary)
  \(Y\). We use \(\Gamma (Y;V)=\Gamma (V)\) to denote the space of smooth
  sections of \(V\), and use \(C^\infty_0(Y; V)\) to denote the space
  of smooth sections whose support lie in compact subspaces in the
  interior of \(Y\).
  Given \(\mathfrc{s}\in
  \Gamma (Y; V)\),
  \(\ds \|\mathfrc{s}\|_{L^p(Y;
    V)}:=\Big(\int_Y|\mathfrc{s}|^p\Big)^{1/p}\). This is sometimes abbreviated as
  \(\|\mathfrc{s}\|_{L^p}\) or \(\|\cdot\mathfrc{s}\|_{p}\). Given a
  connection \(A\) on \(V\), 
  \(\ds \|\mathfrc{s}\|_{L^p_{k/A}(Y;
    V)}:=\Big(\sum_{i=0}^k\int_Y|\nabla_A^k\mathfrc{s}|^p\Big)^{1/p}\),
  where \(\nabla_A\) is used to denote covariant derivatives with
  respect to connections induced from \(A\) and the Levi-Civita
  connection on \(T^*Y\). It is also abbreviated as
  \(\|\frcs\|_{L^p_{k/A}}\) or \(\|\frcs\|_{p,k/A}\). The connection
  \(A\) is sometimes omitted from the notation when its choice is
  obvious or insignificant.  For example, when \(V=\ul{\bbR}\) is the trivial
  \(\bbR\)-bundle, then \(L^p_k(Y)\) denotes \(L^p_{k/A}(Y;
  \ul{\bbR})\) when \(A\) is taken to be the trivial connection.

  Let \(L^p_{k/A, loc}(Y;
    V)\) denote the space consisting all sections of \(V\) whose
    restriction to any compact subspace of \(Y\) is in \(L^p_{k/A}\).
\item Given a topological space \(M\), \(\mathring{M}\) denotes the
  interior of \(M\).
  \item \(C_*, C'_*\) with various subscripts \(*\) usually denote
    positive constants whose precise values are not important, and
    possibly vary with each occurrence. Similarly for \(r_0\). 
    
   \end{itemize}

   This article frequently refers to various literature, which
unfortunately use different conventions. For the reader's convenience,
we clarify some of their relations here. The Seiberg-Witten equations
in this article follow the convention of \cite{KM} and \cite{L3}. In Taubes's articles, \(F_A/2\) above is
replaced by \(F_A\). This results in a difference of factor 2 in many
expressions below from their analogs in Taubes's articles.
To sum up, 
\[\begin{split}\Psi & =\Psi_{\rm KM}=\Psi_{\rm LT}/\sqrt{2}=\psi _{\rm L1}/\sqrt{2}\\
\frac{i\mu}{4} &=irw_{f}|_{\rm LT}=-2w|_{\rm KM}=-\frac{i}{2}\omega|_{\rm L1};\\
[\rho ^{-1}(\Psi ^*\Psi )_0]& =[\rho ^{-1}(\Psi ^*\Psi
)_0]_{\rm KM}=-[\Psi ^\dag\tau \Psi ]_{\rm LT}
\end{split}
\]
where the first expressions in all three lines are in the notation used in
this article, and the subscripts \(KM\), \(PFH\), \(har\) refer
respectively to their counterparts in \cite{KM}, \cite{LT}, and
\cite{L1}. 

\paragraph{Acknowledgement} This work is supported in part by Hong
Kong RGC grant GRF-14301622.

\section{Preliminaries}

\subsection{Some definitions}

Let \(\chi (t)\) denote a non-negative, non-decreasing smooth real function
  on \(\bbR\) such that \(\chi (t)=0\) on \((-\infty, 1]\), and \(\chi
  (t)=1\) on \([2, \infty)\). Let \(\bar{\chi }(t):=\chi (t)+\chi
  (-t)\). Let \(\chi _{e}\co Y\to \bbR\) be the nonnegative function
  defined by
  \[
    \chi_{e}=\begin{cases}
      \iota^*\pi _\bbR^*\bar{\chi } & \text{on \(Y\backslash Y_c\),
        where \(\pi _\bbR\co \bbR^\pm_t\times \Sigma _\pm\to\bbR^\pm_t\)
        denotes the projection,
        }\\
      0 & \text{on \(Y_c\),}
    \end{cases}
  \]
  and let \(\tilde{t}\co Y\to \bbR\) denote the function defined by
  \[
    d\tilde{t}=\begin{cases}
     \chi _{e}\, \iota^*\pi _\bbR^*dt & \text{on \(Y\backslash Y_c\), 
        }\\
      0 & \text{on \(Y_c\),}
    \end{cases}; \quad \tilde{t}\Big|_{Y_c}=0.
    \]
\begin{defn}[Weighted Sobolev norms]\label{def:ext-sobolev}
  Let \(V\to Y\) be a euclidean or hermitian vector bundle over a CCE
  \(Y\). Fix \(\epsilon \in \bbR\). Given \(\mathfrc{s}\in \Gamma (Y; V)\),
  \(\ds \|\mathfrc{s}\|_{L^p_{:\epsilon }(Y;
    V)}:=\Big(\int_Y\Big|e^{\epsilon|\tilde{t}|}\mathfrc{s}\Big|^p\Big)^{1/p}\). This is sometimes abbreviated as
  \(\|\mathfrc{s}\|_{L^p_{:\epsilon }}\) or
  \(\|\cdot\mathfrc{s}\|_{p:\epsilon }\). Given a euclidean/hermitian 
  connection \(A\) on \(V\), 
  \(\ds \|\mathfrc{s}\|_{L^p_{k/A:\epsilon }(Y;
    V)}:=\Big(\sum_{i=0}^k\int_Y\Big|e^{\epsilon|\tilde{t}|}\nabla_A^k\mathfrc{s}\Big|^p\Big)^{1/p}\),
  where \(\nabla_A\) is used to denote covariant derivatives with
  respect to connections induced from \(A\) and the Levi-Civita
  connection on \(T^*Y\). It is also abbreviated as
  \(\|\frcs\|_{L^p_{k/A:\epsilon }}\) or \(\|\frcs\|_{p,k/A:\epsilon }\). The connection
  \(A\) is sometimes omitted from the notation when its choice is
  obvious or insignificant.
   \(L^p_{:\epsilon }(Y;
    V)\) and \(L^p_{k/A:\epsilon }(Y;V)\) denote respectively the
    Banach spaces resulting from completing \(C^\infty_0(Y; V)\) with
    respect to the norms \(\|\cdot\|_{p:\epsilon }\) and
    \(\|\cdot\|_{p,k/A:\epsilon }\).
    \item \(w\) is in the cohomology class \(2\pi c_1(\grs)\).
\end{defn}

\begin{defn}[Extended (weighted) Sobolev space]\label{def:ext-sobolev}
Let \(Y\co \Sigma _-\to \Sigma _+\) be a CCE.  Let \(\pi _2\) denote the projection to the second factor of the
  product \((-\infty, 0)\times \Sigma _-\) or \((0, \infty ) \times
  \Sigma _+\), and let \(\pi _\Sigma \co \scrE^0_\pm[Y]\to \Sigma
  _\pm\) be given by \(\pi _\Sigma :=\pi _2\circ\iota^{-1}\), where \[\text{
  \(\scrE^R_-[Y]=\iota\big((-\infty, -R)\times \Sigma _-\big)\),
  \(\scrE^R_+[Y]=\iota\big((R,\infty)\times \Sigma _+\big)\).}\]

  Let \(V\to  Y\) be a euclidean/hermitian vector bundle with
  bundle isomorphisms \(\iota_V\co \pi _2^*V_\pm\to
  V\Big|_{\scrE'_\pm[Y]}\), where \(V_\pm\to\Sigma _\pm\) are
  euclidean/hermitian vector bundles.
  \[\begin{CD}
      \pi _2^*V_\pm   @>\iota_V>> V\Big|_{\scrE^0_\pm[Y]}\\
      @VVV @VVV\\
      \bbR^\pm\times\Sigma _\pm  @>\iota>> \scrE^0_\pm[Y]
    \end{CD}
  \]
  Fix eucliean/hermitian connections
  \(A_\pm\) on \(V_\pm\). Let \(A_0\) be a euclidean/hermitian
  connection on \(V\) such that it agrees on the induced connection
  from \(A_\pm\) over \( V\Big|_{\scrE'_\pm[Y]}\). Let \(A\) be a  euclidean/hermitian
  connection on \(V\) such that \(A-A_0\in L^p_{l}\), \(l>3/p\),
  \(l\geq k\).  Let \(\chi _{e\pm}\co Y\to \bbR\) be
  \[
    \chi_{e\pm}:=\begin{cases}
      \chi _e & \text{on \(\scrE^0_\pm[Y]\),}\\
      0 & \text{on \(Y\backslash \scrE^0_\pm[Y]\).}
    \end{cases}
  \]
  Given \((\frcs_-, \frcs_+)\in L^p_{k/A_-}(\Sigma _-; V_-)\times
  L^p_{k/A_+}(\Sigma _+; V_+)\), let \( \bfe_{(\frcs_-, \frcs_+)}\in
  L^p_{k/A, loc}(Y; V)\) be given by
  \[
    \bfe_{(\frcs_-, \frcs_+)}:=\chi _{e+}\pi _\Sigma ^*\frcs_++\chi _{e-}\pi _\Sigma
  ^*\frcs_-. 
    \]
  Then given \(\epsilon \geq 0\), let \( \hat{L}^p_{k/A
    :\epsilon }(Y; V)\) denote the space 
  \[
\hat{L}^p_{k/A
    :\epsilon }(Y; V):=  \{\frcs \, |\, \text{
  \(  \exists\frcs_\pm\in L^p_{k/A_\pm }(\Sigma _\pm; V_\pm)\)
  s.t. \(\frcs-\bfe_{(\frcs_-, \frcs_+)}\in L^p_{k/A :\epsilon }(Y; V)\).}\}
    \]
  Let \(\Pi _{\pm\infty}\co  \hat{L}^p_{k/A
    :\epsilon }(Y; V)\to L^p_{k/A_\pm }(\Sigma _\pm; V_\pm)\) denote
  the epimorphism given by \(\frcs \mapsto \frcs_\pm\). Given a
  subspace 
  \(\bbL\subset L^p_{k/A_-}(\Sigma _-; V_-)\times L^p_{k/A_+}(\Sigma _+; V_+)\), let \[\hat{L}^p_{k/A
    :\epsilon }(Y; V|\, \bbL):=\big(\Pi _{-\infty}\times \Pi
  _{+\infty}\big)^{-1}\bbL.\]
By construction, \(\bbL\) is a Banach subspace of \(L^p_{k/A_-}(\Sigma
_-; V_-)\times L^p_{k/A_+}(\Sigma _+; V_+)\), and 
\(\big(\Pi _{-\infty}\times \Pi _{+\infty}\big)\co \hat{L}^p_{k/A
    :\epsilon }(Y; V|\, \bbL)\to \bbL\) is a Banach  bundle over
  \(\bbL\), with fibers isomorphic to the Banach space \(L^p_{k/A
    :\epsilon }(Y; V)\).  We endow \(\hat{L}^p_{k/A
    :\epsilon }(Y; V|\, \bbL)\) with the topology induced from the
  Banach topology on its fibers and base. 

In the case when \(V=T^*Y\), we identify
\(T^*Y\Big|_{\scrE^0_\pm[Y]}\) with \(\pi _2^*V_\pm\), where
\(V_\pm=\ul{\bbR}\oplus T^^*\Sigma _\pm\), where \(\ul{\bbR}\) is the
trivial bundle spanned by \(dt\). In this manner, we regard \(L^p_k
(\Sigma _\pm; T^*\Sigma _\pm)\) as a subspace of \(L^p_k(\Sigma
_\pm; V_\pm)\) by regarding \(T^*\Sigma _\pm\) as a subbundle of
\(V_\pm\). We define
\[
  \tilde{L}^p_{k:\epsilon }(Y; T^*Y):=\hat{L}^p_{k:\epsilon }\big(Y; T^*Y |L^p_k
(\Sigma _-; T^*\Sigma _-)\times L^p_k
(\Sigma _+; T^*\Sigma _+)\big).
  \]
Let \( \tilde{L}^p_{k:\epsilon }(Y; \bigwedge^2T^*Y)\) be similarly
defined: This time identify \(V_\pm \) with \(T^*\Sigma _\pm\oplus
\bigwedge^2T^*\Sigma _\pm\), and use this splitting to identify
\(L^2_{l}(\Sigma _\pm; \bigwedge^2T^*\Sigma _\pm)\) as a subspace of \(L^2_{l}(\Sigma _\pm; V_\pm)\).
\end{defn}

 \begin{defn}\label{def:adm-w}
    Fix a \(\Spin^c\) CCE \((Y, \grs)\). Given \(l\in \bbN\), an \(l\)-{\em admissible 2-form} \(w\in \Omega^2(Y)\) is a closed 2-form
    satisfying the following conditions:
    \begin{itemize}
      \item \(w\in \tilde{L}^2_{l:\epsilon }(Y; \bigwedge^2T^*Y)\), where
        \(\epsilon >0\) satisfies (\ref{def:e});
      \item \(w\) is in the cohomology class \(4\pi c_1(\grs)\).
        \end{itemize}
        \(w\) is said to be {\em admissible} if it is \(l\)-admissible
        for all \(l\in \bbN\).
      \end{defn}
      Let \(\hat{\scrW}_l=\hat{\scrW}_{l, \grs}\) denote the space of
      \(l\)-admissible 2-forms, and  let
      \(\hat{\scrW}=\bigcap_l\hat{\scrW}_l\). Let \(\scrW^\pm\)/\(\scrW^\pm_l\) denote 
      the space of smooth/\(L^2_l\) closed 2-forms on \(\Sigma _\pm \) in the
      cohomology class \(4\pi c_1^\pm\), where \((c_1^-, c_1^+)\in
      H^2(\Sigma _-;\bbZ)\oplus H^2(\Sigma _+;\bbZ)=H^2(\partial
      Y_c;\bbZ)\) is the image of \(c_1(\grs)\) under the pullback map
      \(\imath_c^*\co H^2(Y;\bbZ)\to H^2(\partial Y_c; \bbZ)\), where
      \(\imath_c\co \partial Y_c\to Y\) is the embedding. The
      end-point maps \(\Pi _{+\infty}\) define a bundle structure on
      \(\hat{\scrW}\) and \(\hat{\scrW}_l\):
\[\big(\Pi _{-\infty}\times \Pi _{+\infty}\big)\co
  \hat{\scrW}\to \scrW^-\times \scrW^+; \quad \big(\Pi _{-\infty}\times \Pi _{+\infty}\big)\co
  \hat{\scrW}_l\to \scrW^-_l\times \scrW^+_l. \]
By construction, the fibers of \(\hat{\scrW}_l\) are 
      affine spaces under the space of exact
      \(L^2_{l:\epsilon }\) 2-forms on \(Y\), denoted as
      \(\rmW_l\). We endow \(\hat{\scrW}_l\) with the topology induced
      from the Banach topologies on its base and fibers, 
and similarly endow
      \(\hat{\scrW}\) with topology induced from the Fr\'echet
      topologies of \(\scrW^\pm\) and \(\rmW:=\bigcap_l\rmW_l\).
 
      \begin{defn}\label{def:adm-config}
    Fix a \(\Spin^c\) CCE \((Y, \grs)\), and let \(\bbS\) denote the
    associated spinor bundle. Fix an isomorphism \(\iota_\bbS\co  \pi
    _2^*\bbS_{\Sigma _\pm   }\to \scrE^0_\pm[Y]\)
\[\begin{CD}
      \pi _2^*\bbS_{\Sigma _\pm   }@>\iota_\bbS>> \bbS\Big|_{\scrE^0_\pm[Y]}\\
      @VVV @VVV\\
      \bbR^\pm\times\Sigma _\pm  @>\iota>> \scrE^0_\pm[Y].
    \end{CD}
  \]
    Choose a reference connection $A_0\in
    \Conn (\bbS)$ such that its restriction to \(\scrE^0_\pm[Y]\)
    agrees with a pull-back connection  \(\pi _2^*B_{0, \pm}\),
    \(B_{0, \pm}\in \op{Conn} (\bbS_{\Sigma _\pm   })\). Given \(l\in
    \bbN\), we say that \(A\in \Conn (\bbS)\) is {\em \(l\)-admissible}
      if \(A_t-(A_0)_t\in \tilde{L}^2_{l}(Y, iT^*Y)\). 
    
Let \(l\in \bbN\), \(l\geq 2\). A configuration $(A, \Psi )$ is {\em \(l\)-admissible} if it satisfies: 
\begin{enumerate}
\item \(A\) is \(l\)-admissible and $\Psi \in \hat{L}^2_{l/A_0}(Y;
  \bbS)$. Note that by Sobolev embedding, when \(A\) is
  \(l\)-admissible with \(l\geq 2\), $\Psi \in \hat{L}^2_{l/A_0}(Y;
  \bbS)$ iff $\Psi \in \hat{L}^2_{l/A}(Y;\bbS)$. 
\item $\big(\frac{\rho(d\ff)}{|d\ff|}-i\big)\Psi \Big|_{\scrE_\pm [Y]}\in
  L^2_{l /A} (\scrE_\pm [Y]; \bbS)$. Note that
  \(\frac{\rho(d\ff)}{|d\ff|}\) is well-defined on the ends of \(Y\)
  as the zero locus of \(d\ff\) falls in \(Y_c\).
\end{enumerate}
$(A, \Psi )$ {\em admissible} if it is admissible \(\forall l\).

By Condition 1 above, there are end-point maps \(\Pi _{\pm\infty}\) from the space of
 \(l\)-admissible configurations to \(\Conn_l   (\bbS_{\Sigma _\pm
 })\times L^2_{l/B_{0, \pm}}(\Sigma _\pm; \bbS_{\Sigma _\pm
 })\), where  \[
   \Conn_l   (\bbS_{\Sigma _\pm
 }):=\{B_{0, \pm}+b|\, b\in L^2_l(\Sigma _\pm ; iT^*\Sigma
 _\pm)\}.\]
Meanwhile, use \(\rho (dt)\) to split
 \(\bbS\Big|_{\scrE^0_\pm[Y]}=\hat{E}\oplus \hat{E}'\), where
 \(\hat{E}\) is the eigenbundle of \(\rho (t)\) with eigenvalue
 \(-i\).  This induces a splitting of
 \begin{equation}\label{split:S-end}
   \bbS_{\Sigma _\pm}=E_{\Sigma
   _\pm}\oplus  E_{\Sigma
   _\pm}\otimes T^{1,0}\Sigma _\pm
\end{equation}
via the bundle isomorphism
 \[\iota_\bbS\co \pi _2^*\bbS_{\Sigma _\pm   }=\pi _2^*E_{\Sigma
   _\pm}\oplus  \pi _2^*E_{\Sigma
   _\pm}\otimes T^{1,0}\Sigma _\pm \to
 \bbS\Big|_{\scrE^0_\pm[Y]}=\hat{E}\oplus \hat{E}'.\]
Condition 2 above implies that \(\Pi _{\pm\infty}\) maps an admissible
to an element \((B_\pm, \Phi _\pm)\), where the \(E_{\Sigma
   _\pm}\otimes T^{1,0}\Sigma _\pm\)-component of \(\Phi _\pm\) under
 the splitting (\ref{split:S-end}) vanishes. Thus, we may identify
 \(\Phi _\pm\) as a section of \(E_{\Sigma_\pm}\), and take the
 codomain of \(\Pi _{\pm\infty}\) to be \(\Conn_l   (\bbS_{\Sigma _\pm
 })\times L^2_{l/B^E_{0, \pm}}(\Sigma _\pm; E_{\Sigma _\pm})\), where
 \(B^E_{0, \pm}:=B^E_{0, \pm}\Big|_{E_{\Sigma _\pm}}\).
 
 Define the  end-point maps \(\Pi _{\pm\infty}\) from the space of
 admissible configurations similarly. 
\end{defn}

\subsection{Existence and genericity of admissible functions}

\begin{prop}\label{prop:adm-f}
  Let \(Y\co \Sigma _-\to \Sigma _+\) be a CCE. Then there exists a
  harmonic function \(\ff\) on \(Y\) satisfying Condition (2) of Definition \ref{def:adm-f}.
  Moreover, any two such functions
  differ by a constant function. Given \(\epsilon >0\) satisfying
  (\ref{def:e}) and a non-negative integer \(k\), there are 
  constants \(C_{\ff} >0\), \(\rmf_\pm\) (depending on \(\ff\)) such
  that the following pointwise bound holds:
  \begin{equation}\label{decay-f}
    \sum_{i=1}^k|\nabla^k(\ff-\tilde{t})\, |+|\ff-\tilde{t}-\rmf_\pm
    |\leq C_{\ff} \,\,  e^{-\epsilon
      |\tilde{t}|} \quad \text{over \(\scrE_\pm [Y]\).}
  \end{equation}
\end{prop}
\pf Consider the differential operator \(D\co \Omega ^0(Y)\oplus\Omega
^1(Y)\to  \Omega ^0(Y)\oplus\Omega
^1(Y)\) given by
\[
  D:=\left[\begin{array}{cc}
       0 & d^*\\
     d & *d\end{array}\right].
\]
Then \(D\) is formally \(L^2\) self-adjoint, and \(D (f, \theta )=0\)
implies that both \(f\) and \(\theta \) are harmonic. \(D\) is of the
Atiyah-Patodi-Singer (APS \cite{APS}) type: Over
\(\iota^{-1}\scrE^0_\pm [Y]=\bbR^\pm_t\times \Sigma _\pm\),
identify each element in \(\Omega ^0(\bbR^\pm_t\times \Sigma _\pm)\oplus
\Omega ^1(\bbR^\pm_t\times \Sigma _\pm)\) with a family of elements in
\(\Omega ^0(\Sigma _\pm)\oplus\Omega ^0(\Sigma _\pm)\oplus\Omega
^0(\Sigma _\pm)\) parametrized by \(t\in \bbR^\pm\) as follows: Assign
to each \[\big(f(t, z), \theta (t, z)=\vartheta(t, z) dt+\theta _z(t, z)\big)\in  \Omega ^0(\bbR^\pm_t\times \Sigma _\pm)\oplus
  \Omega ^1(\bbR^\pm_t\times \Sigma _\pm)\]
the family \[t\mapsto
(f(t, \cdot), \vartheta(t,\cdot), \theta _z(t, \cdot))\in \Omega ^0(\Sigma _\pm)\oplus\Omega ^0(\Sigma _\pm)\oplus\Omega
^1(\Sigma _\pm)=\Gamma (\Sigma _\pm; \ul{\bbR}\oplus \ul{\bbR}\oplus T^*\Sigma _\pm).\]
In the above, \(t\in \bbR^\pm\), \(z\in \Sigma _\pm\), \(\vartheta
(t,\cdot)\in \Omega ^0(\Sigma _\pm)\), and \(\theta _z(t, \cdot)\in
\Omega ^1(\Sigma _\pm)\). 
Then under the aforementioned  identification, 
\[\iota^*\circ D\circ (\iota^*)^{-1}=\sigma  \big(\frac{d}{dt}+B\big),\]
where \(\sigma \co E_\pm \to E_\pm\) is a bundle automorphism, \(E_\pm := \ul{\bbR}\oplus
\ul{\bbR}\oplus T^*\Sigma _\pm\):
\[
  \sigma =\left[\begin{array}{ccc} 0 &-1 &0\\
                  1 & 0 & 0\\
                  0 & 0 & *_z
                \end{array}
              \right]; \]
and \(B\co \Gamma (\Sigma _\pm ; E_\pm)\to \Gamma (\Sigma _\pm ; E_\pm)\) is the differential
operator            
            \[
B_\pm= \left[\begin{array}{ccc} 0 & 0 &-*_z d_z\\
                  0 & 0 & d_z^*\\
                   *_z d_z & d_z & 0
                \end{array}
              \right].            \]
In the above, \(*_z\co \Omega ^*(\Sigma _\pm)\to \Omega ^{2-*}(\Sigma
_\pm)\) denote the 2-dimensional Hodge dual; \(d_z=d\co  \Omega ^*(\Sigma _\pm)\to \Omega ^{*}(\Sigma
_\pm)\) denotes the 2-dimensional exterior derivative.

Note that \(\sigma \) and \(B_\pm\) satisfy the properties that
\[\text{
\(\sigma ^2=-1\), \(\sigma ^*=-\sigma \), \(\sigma
B_\pm+B_\pm\sigma =0\), }\]
and \(B\) is formally \(L^2\) self-dual adjoint.  It extends
to a self-adjoint Fredholm operator denoted by the same notation: 
\[
  B_\pm\co L^2_1(\Sigma _\pm; E_\pm)\to  L^2(\Sigma _\pm, E_\pm). 
\]
The kernel and the cokernel of \(B_\pm\) are both
\[
  \bbH_{B_\pm}=\{ (f, \vartheta, \theta _z) \, |\, f, \vartheta,
  \theta _z \, \text{harmonic}\}\simeq H^0(\Sigma _\pm)\oplus
  H^0(\Sigma _\pm)\oplus H^1(\Sigma _\pm). 
\]
Moreover, as observed in \cite{CLM}, \(\sigma \) induces a symplectic
form \(\Omega ^\sigma _\pm\) on \(L^2(\Sigma _\pm, E_\pm)\):
\[
  \Omega ^\sigma _\pm (h, h'):=\langle h, \sigma  h'\rangle_{L^2}, 
\]
which restricts to a symplectic form (denoted by the same notation) on
\(\bbH_{B_\pm}\). Let \(\bar{\Omega }^\sigma \) denote the symplectic
form  \( (-\Omega ^\sigma _-)\oplus \Omega ^\sigma _+\) on
\(L^2(\Sigma _-, E_-)\oplus L^2(\Sigma _+, E_+)\), which in turn
induces a symplectic form on \(\bbH_{B_-}\times \bbH_{B_+}\) denoted
by the same notation. Let \(\bar{\sigma }:=(-\sigma)\oplus
\sigma\). Then \(\bar{\sigma }\) defines a complex structure on
\(L^2(\Sigma _-, E_-)\oplus L^2(\Sigma _+, E_+)\) compatible with \(\bar{\Omega }^\sigma \), which in turn
induces a complex structure on \(\bbH_{B_-}\times \bbH_{B_+}\)
compatible with \(\bar{\Omega }^\sigma \), also denoted
by the same notation. 

As a self-adjoint operator, \(B_\pm\) has a discrete spectrum
\(\op{Spec}(B_\pm)\) in the real
line. Fix
\begin{equation}\label{def:e}
  \text{ \(\epsilon >0\) such that \(\epsilon  <\min \Big(\min_{\lambda \in
  \op{Spec}(B_+), \lambda \neq 0},  |\lambda  | \min_{\lambda \in
  \op{Spec}(B_-), \lambda \neq 0} |\lambda  |\Big)=:\epsilon _0\).}
\end{equation}
Let \[
  D_{:\epsilon }\co L^2_{1 :\epsilon } (Y;
\ul{\bbR}\oplus  T^*Y)\to   L_{:\epsilon }^2(Y;
\ul{\bbR}\oplus  T^*Y)
\]
denote the operator obtained by completing  \(D\co C^\infty_0(Y;
\ul{\bbR}\oplus  T^*Y)\to  C^\infty_0(Y;
\ul{\bbR}\oplus  T^*Y)\) with respect to the \(L^2_{1:\epsilon
}\)-norm.

Note that for any  \(\bar{h}=(h_-, h_+)\in  \bbH_{B_-}\times
\bbH_{B_+}\),   \(D\bfe_{\bar{h}}\) is compactly supported; so
\(D_{:\epsilon }\) extends to define an operator
\[
  \hat{D}_{:\epsilon }\co \hat{L}^2_{1:\epsilon }(Y;
\ul{\bbR}\oplus  T^*Y |\, \bbH_{B_-}\times \bbH_{B_+})\to L^2_{:\epsilon }(Y;
\ul{\bbR}\oplus  T^*Y).
\]
Given a subspace \(L\subset \bbH_{B_-}\times
\bbH_{B_+}\), let \(\hat{D}_{:\epsilon  |L}\co \hat{L}^2_{1:\epsilon }(Y;
\ul{\bbR}\oplus  T^*Y |\, L)\to L^2_{:\epsilon }(Y;
\ul{\bbR}\oplus  T^*Y)\) denote the restriction of \(
\hat{D}_{:\epsilon }\) to \(\hat{L}^2_{1:\epsilon }(Y;
\ul{\bbR}\oplus  T^*Y |\, L) \subset \hat{L}^2_{1:\epsilon }(Y;
\ul{\bbR}\oplus  T^*Y |\, \bbH_{B_-}\times \bbH_{B_+})\). 


\begin{lemma}\label{lem:fred-D}
  Fix \(\epsilon \) satisfying (\ref{def:e}). Then:
  
\noindent {\rm (1)} The operators  \(D_{ :\epsilon }\), \(\hat{D}_{:\epsilon  |L}\) are
Fredholm, where \(L\) is an arbitrary subspace of \(\bbH_{B_-}\times
\bbH_{B_+}\).

\noindent {\rm (2)} Let \(\bbH_Y:=\ker \hat{D}_{:\epsilon  }\), then
\[
  \bbH_Y=\{(C, \theta _h)\, |\, \text{\(C\) is a constant function;
    \(\theta _h\) is a harmonic 1-form on \(Y\)}\}.\]
Moreover, \(L_D:=(\Pi
_{-\infty}\times\Pi _{+\infty})\bbH_Y\) is a
Lagrangian subspace in \((\bbH_{B_-}\times \bbH_{B_+}, \bar{\Omega
}^\sigma )\).  The fiber of the surjection \(\Pi
_{-\infty}\times\Pi _{+\infty}\co \bbH_Y\to L_D\) is isomorphic to
the image of \(H^1(Y_c, \partial Y_c)\) in \(H^1(Y_c)\).

\noindent {\rm (3)} \(\hat{D}_{:\epsilon  |\bar{\sigma }L_D}\) is of
index 0, whose kernel and cokernel are both isomorphic to the image of \(H^1(Y_c, \partial Y_c)\) in \(H^1(Y_c)\).
\end{lemma}

\pf
(1) Since \(\hat{D}_{:\epsilon  |L}\) is a finite dimensional
extension of \(D_{:\epsilon }\) for any \(L\), it suffices to verify
that \(D_{:\epsilon }\) is Fredholm. This follows from the argument in
\cite{APS}, noting that when \(\epsilon \) satisfies (\ref{def:e}),
the parametrix \(R\) constructed in p.54 of \cite{APS} is also a
parametrix for \(D_{:\epsilon }\).

(2) The first statement follows from Proposition 3.15 of \cite{APS},
noting that \(D^2=d^*d+dd^*\), and the observation that when \(\epsilon \) satisfies
(\ref{def:e}), an extended \(L^2\)-solution \(\frcs\) of \(D\frcs =0\) in the
sense of \cite{APS} is in \(\hat{L}^2_{1:\epsilon }\), since over
\(\scrE^0_\pm[Y]\), \(\frcs\) takes the form 
\[\frcs \Big|_{\scrE^0_\pm[Y]}=\iota^* \sum_{\lambda \in \op{Spec}
    (B_\pm), \pm\lambda \geq 0}e^{-\lambda t}\xi_\lambda , \]
where \(\xi_\lambda \) is an eigenfunction of \(B_\pm\) with
eigenvalue \(\lambda \).

The second statement follows from \cite{CLM} Proposition 2.3.

The third statement follows from the following observations: Each
fiber of \(\Pi _{-\infty}\times \Pi _{+\infty}\)
is an affinement space over \[
  \ker D_{:\epsilon }=\{(0, \theta _h) |\, \text{ \(\theta_h\in
    L^2_{:\epsilon}\) is a harmonic 1-form}\},\]
and the space of \(L^2_{:\epsilon}\) harmonic 1-forms on \(Y\) agrees
with the space of \(L^2\) harmonic 1-forms, since both \(\ker
D_{:\epsilon }\) and \(\ker D_{:\epsilon 0}\)
consist of elements \(\frcs\) with taking the form
\[\frcs \Big|_{\scrE^0_\pm[Y]}=\iota^* \sum_{\lambda \in \op{Spec}
    (B_\pm), \pm\lambda >0}e^{-\lambda t}\xi_\lambda , \]
where \(\xi_\lambda \) is again an eigenfunction of \(B_\pm\) with
eigenvalue \(\lambda \). We denote this space as \(\scrH^1_c(Y)\). Finally, the latter space is isomorphic to
the image of \(H^1(Y_c, \partial Y_c)\) in \(H^1(Y_c)\) according to
\cite{APS} Proposition 4.9.

(3) Observe that \[\begin{split}
  \ker \hat{D}_{:\epsilon  | \bar{\sigma }L_D}& =\bbH_Y\cap(\Pi
  _{-\infty}\times \Pi _{+\infty})^{-1}\bar{\sigma }L_D\\
  & = \bbH_Y\cap(\Pi
  _{-\infty}\times \Pi _{+\infty})^{-1}(\bar{\sigma }L_D\cap L_D)
  \\& = \bbH_Y\cap(\Pi
  _{-\infty}\times \Pi _{+\infty})^{-1}(0)\\
  &=\ker D_{:\epsilon }=\{0\}\oplus\scrH^1_c(Y). 
\end{split}
\]
Regard \(L^2_{:\epsilon }\) as a Hilbert space with inner product
\[
  \langle f, g\rangle_{2:\epsilon }:=\langle e^{\epsilon |\tilde{t}|}
  f, e^{\epsilon |\tilde{t}|} g\rangle_2, 
\]
where
    \(\langle\cdot, \cdot\rangle_2=\langle\cdot, \cdot\rangle_{L^2}\)
    denotes the \(L^2\) inner product. Then \(q\in \op{coker}
    D_{:\epsilon }\) iff
    \[
      \langle D\frcs, q\rangle_{2:\epsilon }=\langle  \frcs,
      D(e^{2\epsilon |\tilde{t}|} q)\rangle_2=0 \quad \forall \frcs\in L^2_{1:\epsilon }.
    \]
    Since \(C^\infty_0\) is dense in both \(L^2\) and
    \(L^2_{1:\epsilon }\), this implies that \(e^{2\epsilon
      |\tilde{t}|} q\in L^2_{:-\epsilon }\) is harmonic. The argument
    in part (2) above implies that such an element is in
    \(\hat{L}^2_{1:\epsilon }\), and hence \begin{equation}\label{eq:cok-D}
      \op{coker}
    D_{:\epsilon }=\{ e^{-2\epsilon |\tilde{t}|} h|\, h\in
    \bbH_Y\}.
  \end{equation}
  We claim that \begin{equation}\label{eq:cok-hatD}
    \op{coker} \hat{D}_{:\epsilon  |
      \bar{\sigma }L_D}=\{ e^{-2\epsilon |\tilde{t}|} (0, \theta _h)|\, \theta _h\in
    \scrH_c^1(Y)\}\simeq \scrH_c^1(Y).
  \end{equation}
  
    Since \(\hat{L}^2_{1:\epsilon }(Y; \ul{\bbR}\oplus  T^*Y |\,
    \bar{\sigma }L_D)=\op{Span} \{\bfe_{\bar{h}}|\, \bar{h}\in
      \bar{\sigma }L_D\}\oplus L^2_{1:\epsilon }(Y; \ul{\bbR}\oplus
      T^*Y)\), it suffices to show that:
\begin{itemize}
  \item[(i)]     For each \(\bar{h}\in L_D\), \(\bar{h}\neq 0\), 
      there exists \(h\in \bbH_Y\) such that \(\langle
      D\bfe_{\bar{\sigma }\bar{h}}, e^{-2\epsilon |\tilde{t}|} h\rangle_{2:\epsilon
      }\neq 0\);
  \item[(ii)]  \(\langle D\bfe_{\bar{\sigma }\bar{h}}, e^{-2\epsilon
      |\tilde{t}|} (0, \theta _h)\rangle_{2:\epsilon
      }=0\) \(\forall\)  \(\bar{h}\in L_D\),  \(\theta _h\in
      \scrH_c^1(Y)\). 
    \end{itemize}
    Both of the statements above follows from the following
    computation:  
  Given \(\rmh\in \bbH_Y\), let \(\rmh_{\pm\infty}:=\Pi
  _{\pm\infty}\rmh\). Write \(\bar{h}=(h_-, h_+)\). By the Stokes'
  theorem, 
  \[\begin{split}
    \langle D\bfe_{\bar{\sigma }\bar{h}} , & e^{-2\epsilon |\tilde{t}|}
    \rmh\rangle_{2:\epsilon} \\
    &=\langle D\bfe_{\bar{\sigma }\bar{h}}, \rmh\rangle_{2}\\
    & =\langle      \bfe_{\bar{\sigma }\bar{h}}, Dh\rangle_{2}+\langle \sigma \sigma  h_+,
      \rmh_{+\infty}\rangle_{L^2(\Sigma _+; E_+)}-\langle \sigma (-\sigma  h_-),
      \rmh_{-\infty}\rangle_{L^2(\Sigma _-; E_-)}\\
      &=-\langle h_+, \rmh_{+\infty}\rangle_{L^2(\Sigma _+; E_+)}-\langle h_-,
      \rmh_{-\infty}\rangle_{L^2(\Sigma _-; E_-)}. 
    \end{split}\]
  To verify (i), simply take \(\rmh\) to be an element with
  \(\rmh_{\pm\infty}=h_\pm\). To verify (ii), take \(\rmh=(0, \theta
  _h)\), \(\theta _h\in \scrH_c^1(Y)\). Then \(\rmh_{\pm\infty}=0\). 
  \epf

  Return now to the proof of Proposition \ref{prop:adm-f}. We shall
  show that there exists a \(f\in \hat{L}^2_{1:\epsilon }(Y)\), such
  that \(\ff=  \tilde{t}+f\) is a harmonic function satisfying
  Condition (2) of Definition \ref{def:adm-f}. That is, 
  \(d^*df=-d^*d\tilde{t}.\) Note that \(-d^*d\tilde{t}\) is compactly
  supported on \((\scrE^1_-[Y]\backslash\scrE_-[Y])\cup
  (\scrE^1_+[Y]\backslash\scrE_+[Y])\), and \(\displaystyle\int_Y
  (-d^*d\tilde{t})=0\). Thus, by (\ref{eq:cok-D}) and Lemma
  \ref{lem:fred-D} (2), \((-d^*d\tilde{t}, 0)\) is \(L^2_{:\epsilon
  }\)-orthogonal to the cokernel of \(D_{:\epsilon }\); hence there
  exists \((f_0, \theta _0)\in L^2_{1:\epsilon }(Y; \ul{\bbR}\oplus
  T^*Y)\) such that \(D(f_0, \theta _0)=(-d^*d\tilde{t},
  0)\). Morever, the space of all such solutions is an affine space
  under \(\ker D_{:\epsilon }=\{(0, \theta
  _h) |\, \theta _h\in \scrH^1_c(Y)\}\). Thus, we can and shall choose a solution \((f_0,
  \theta _0)\) such that \[
    \langle(f_0, \theta _0), (0, \theta_h)\rangle_2=\langle(f_0, \theta _0), e^{-2\epsilon |\tilde{t}|}(0, \theta
  _h)\rangle_{2:\epsilon }=0 \quad \forall \theta _h\in \scrH^1_c(Y). \]
  Recalling (\ref{eq:cok-hatD}), this implies that \((f_0, \theta
  _0)\) is in the image of \(\hat{D}_{:\epsilon  }\), and thus there
  exists a \((f, \theta )\in  \hat{L}^2_{1:\epsilon }(Y;
  \ul{\bbR}\oplus  T^*Y |\, \bbH_{B_-}\times \bbH_{B_+})\) such that
  \(D(f, \theta )=(f_0, \theta _0)\). Now, \(D^2(f, \theta
  )=(-d^*d\tilde{t}, 0)\) implies that \(\ff :=\tilde{t}+f\) is
  a harmonic function satisfying
  Condition (2) of Definition \ref{def:adm-f}.  Morever, if
  \(\mathpzc{g}\) is another such function,
  then \(\mathpzc{g}-\ff\in \hat{L}^2_1\) and \(d^*d(\mathpzc{g}-\ff)=0\). Thus,
  \(\mathpzc{g}-\ff\) is a constant.

  To verify (\ref{decay-f}), note that since \((-d^*d\tilde{t}, 0)=0\)
  over \(\scrE_\pm [Y]\), \((f, \theta )\) takes the form
  \begin{equation}\label{1}
  (f, \theta )\Big|_{\scrE_\pm[Y]}=\iota^* \Big(\sum_{\lambda \in \op{Spec}
    (B_\pm), \pm\lambda >0}\frac{e^{-\lambda (t\mp 2)}}{\lambda }\xi_\lambda ^\pm+\xi_0^\pm\Big),
\end{equation}
where \(\xi^\pm_\lambda \) is an eigenfunction of \(B_\pm\). Now, \((f,
\theta )\Big|_{\scrE^0_\pm[Y]}-\iota^*\xi_0^\pm\in
L^2_{2:\epsilon}(\scrE^0_\pm[Y]; \ul{\bbR}\oplus  T^*Y)\). On the
other hand,
\begin{equation}\label{}
  \begin{split}
  \Big\|\sum_{\lambda \in \op{Spec}  (B_\pm), \pm\lambda >0}\frac{1}{\lambda }\xi^\pm_\lambda
  \Big\|_{C_k}& \leq \|(f,\theta
  )-\iota^*\xi^\pm_0\|_{C_k(\scrE_\pm^1[Y]\backslash \scrE_\pm^3[Y];
    \ul{\bbR}\oplus  T^*Y)}\\
  & \leq
  C\|(f,\theta)-\iota^*\xi^\pm_0\|_{L^2_{k+2}(\scrE_\pm^1[Y]\backslash
    \scrE_\pm^3[Y]; \ul{\bbR}\oplus  T^*Y)}\\
\end{split}
\end{equation}
by Sobolev embedding, where \(C\) is a (\(\ff\)-independent)
positive constant. By elliptic bootstrapping,
\begin{equation}\label{3}
\begin{split}
  \|(f,\theta)-\iota^*\xi^\pm_0\|_{L^2_{k+2}(\scrE_\pm^1[Y]\backslash
    \scrE_\pm^3[Y]; \ul{\bbR}\oplus  T^*Y)} &\leq C_k \|(f,\theta)-\iota^*\xi^\pm_0\|_{L^2_{2}(\scrE_\pm^1[Y]\backslash
    \scrE_\pm^3[Y]; \ul{\bbR}\oplus  T^*Y)} \\
  & \leq C'_k\|(f,\theta)-\iota^*\xi^\pm_0\|_{L^2_{2:\epsilon
    }(\scrE_\pm^0[Y]; \ul{\bbR}\oplus  T^*Y)}, 
\end{split}
\end{equation}
where \(C_k\) \(C_k'\) are (\(\ff\)-independent) positive
constants. (\ref{decay-f}) now follows from a combination of
(\ref{1})-(\ref{3}) together with the observation that
\(\xi_0^\pm=(\rmf_\pm, \rmh_\pm)\), where \(\rmf_\pm\) are constants,
and \(\rmh_\pm\) are harmonic 1-forms. 
  \epf

  Admissible functions are generic in the following sense.
  \begin{prop}
    Let \(Y\) be a CCE with cylindrical metric \(g_0\), and let
    \(\ff_0\) be a function satisfying Condition (2) of Definition
    \ref{def:adm-f} that is harmonic with respect to \(g_0\). Note that by (\ref{decay-f}), there exists \(R\geq 0\) such that
\(\|d\ff\|_{C_0(\scrE^R_\pm[Y])}> 1/2\). We redefine \(Y_c\) to be
\(Y_c\backslash (\scrE^R_-[Y]\cup \scrE^R_+[Y])\).  Given 
\(\varepsilon>0\), let \[
  \scrU_\varepsilon:=\{h\, |\, h\in
  C^\infty_0(Y_c; \op{Sym}^2T^*Y), \|h\|_{C^2}\leq \varepsilon\},\]
endowed with the Fr\'echet topology. 
Choose \(\varepsilon\) to be  sufficiently small such that \(\forall
h\in \scrU_\varepsilon\), \[
      g_h:=\begin{cases}
        g_0+h &\text{over \(Y_c\)}\\
        g_0 &\text{over \(Y\backslash Y_c\)}
      \end{cases}
    \]
      is also a cylindrical metric on \(Y\). By the previous
      proposition, there exists a 
      unique function \(\ff_h\) satisfying:
      \BTitem \label{adm-morse}
        \item Condition (2) of Definition
          \ref{def:adm-f} holds; 
        \item \(\ff_h\) is harmonic with respect to \(g_h\);
        \item \(\sum_{i=0}^k|\nabla^k(\ff_h-\ff_0)|\leq C_{h} \, e^{-\epsilon
     |\tilde{t}|} \quad \text{over \(\scrE_+[Y]\),}\) where \(C_h>0\) is a
   constant depending on both \(h\) and \(\ff\).  
  \ETitem 
        Then when \(\varepsilon\) is sufficiently small, the zero
        locus of \(d\ff_h\) lies in the interior of \(Y_c\) \(\forall
        h\in \scrU_\varepsilon\), and there exists a Baire subset
        \(\scrU_\varepsilon^{reg}\subset \scrU_\varepsilon\), such
        that \(\ff_h\) is admissible when \(h\in \scrU_\varepsilon^{reg}\).
  \end{prop}
  \pf This follows from a more-or-less standard transversality
  argument via the Sard-Smale theorem. Detailed proofs in similar
  contexts are written down in e.g. \cite{H} (for compact \(Y\)) and
  \cite{Lv2} (for MEE).

  Let \(*_g\) denote the Hodge dual with respect to the metric \(g\),
  and let \(\delta _h*:=*_{g_h}-*_{g_0}\). Then \(f_h:=\ff_h-\ff_0\)
  satisfies:
  \begin{equation}\label{DE:f-pert}
    d*_{g_h}df_h=-d\big((\delta _h*)d\ff_0\big).
  \end{equation}
  Since the integral of the right hand side over \(Y\) equals 0, the
  arguments in the proof of the previous proposition the preceding
  equation has a solution \(f_h\in \hat{L}^2_{2:\epsilon }\), unique
  modulo constant functions. We choose the constant so that the third
  bullet of (\ref{adm-morse}) holds. In \cite{LiT}, a symmetric
  Green's function \(G_g(x,y)\) is constructed for complete Riemannian
  manifolds. This Green's function has the following properties:
  \BTitem\label{def:LiT}
    \item \(G_g(x, y)\sim 4\pi \op{dist} (x, y)^{-1}\) as \(x\to  y\);
     \item \(G_g(x, y)\Big|_{y\in Y\backslash B_x(R)}\) is bounded,
       where \(R>0\) and \(B_x(R)\) is a geodesic ball of radius \(R\)
       centered at \(x\). 
  \ETitem
  Thus, the function
  \[
    \rmf_h(x):=-\int _YG_{g_h}(x,y)d_y\big((\delta _h*_{g_0})d_y\ff_0(y)\big)=\int_Y(d_yG_{g_h}(x,y))(\delta _h*)d_y\ff_0(y)
  \]
  is also a solution to (\ref{DE:f-pert}). (In the above, \(d_x, d_y\)
  respectively denote the exterior derivative in the variable \(x\),
  \(y\).) Moreover,
  \begin{equation}\label{sol:df}
d_x\rmf_h(x)=\int_Y\big(d_xd_yG_{g_h}(x,y)\big)(\delta _h*_{g_0})d_y\ff_0(y)\in
L^2_{:\epsilon }. 
\end{equation}
To see this, note that since when \(x\in \scrE_\pm[Y]\), \(y\in Y_c\), 
\((G_{g_h}(x,y)\Big|_{\scrE^{2|\tilde{t}(x)|}_\pm}, 0)\) is harmonic
and thus by the second bullet
of (\ref{def:LiT}) takes the form \[\sum_{\lambda \in \op{Spec} (B_\pm), \pm\lambda\geq
0}\xi_{\lambda  , y}e^{-\lambda \tilde{t}(x)},\] where for fixed \(y\in Y_c\),
\(\xi_{\lambda , y}\) is an eigenfunction of \(B_\pm\) varying
smoothly with \(y\). Since \(\xi_{0,y}=(C_{0, y}, 0)\), where
\(C_{0,y}\) is a constant function (depending on \(y\)), there is constant
\(C\) such that \(|d_xd_yG_{g_h}(x,y)|\leq Ce^{-\epsilon_0
  |\tilde{t}(x)|}\) \(\forall y\in Y_c\) as \(Y_c\) is
compact. Plugging this into the right hand side of the equation
(\ref{sol:df}), and recalling that \(h\) is compactly supported on
\(Y_c\), we have thus verified that \(d\rmf_h\in L^2_{:\epsilon}\). 

Next, note that both \((0, df_h)\) and \((0, d\rmf_h)\) are
\(L_{2:\epsilon }\) solutions
to \[D_{:\epsilon , g_h}(-)=\Big(-*_{g_h}d\big((\delta _h*_{g_0})d\ff_0\big), 0\Big),\]
and both are \(L^2\)-orthogonal to \(\ker D_{:\epsilon , g_h}\). Thus,
\[df_h=d\rmf_h=\int_Y\big(d_xd_yG_{g_h}(x,y)\big)(\delta
  _h*_{g_0})d_y\ff_0(y).\]
It follows that \(\|df_h\|_{C^0}\leq C\|h\|_{C^0}\) for a positive
constant \(C\). Since \(Y_c\) is chosen such that \(|d\ff_0|>1/2\)
over \(Y\backslash Y_c\), for sufficiently small \(\varepsilon>0\),
\(|d\ff_h|>0\) over  \(Y\backslash Y_c\) \(\forall h\in
\scrU_\varepsilon\). Thus, the zero loci of \(d\ff_h\) lies in the
interior of \(Y_c\) \(\forall h\in
\scrU_\varepsilon\). In particular, since \(Y_c\) is compact, the zero
loci of \(d\ff_h\) is compact and when \(\ff_h\), it consists of
finitely many points.

With the above understood,  a straightforward adaptation of  the argument in Theorem 2.19 in
  \cite{H} shows that there is an open dense subset
  \(\scrU^{reg}_{\varepsilon, l}\) in \[\scrU_{\varepsilon,
    l}:=\{h|\, h\in C^l_0(Y_c; \op{Sym}^2T^*Y), \|h\|_{C^2}\leq
  \varepsilon\}\] for every integer \(l\geq 2\). More explicitly,
modify the argument in \cite{H} as follows: 
\begin{itemize}
\item Replace Equation (32) in \cite{H} with
  \[
    ev_{0, x}(g_h)(\rmh)=\int_Y\big(d_xd_yG_{g_h}(x,y)\big)(\delta
    _{\rmh}*_{g_h})d_y\ff_{h}(y), \]
  where \(x\in  (d\ff_h)^{-1}(0)\subset \mathring{Y}_c\), \(h\in\scrU_{\varepsilon,
    l}\), and \(\rmh\in C^l_0(Y_c; \op{Sym}^2T^*Y)\).   (Note that the \(G_g\) in \cite{H} denotes the
  Green's function for 1-forms instead.)
\item Replace the computation around Equations (33) and (34) of
  \cite{H} by the following. Choose a trivialization of
  \(\bigwedge^2T^*Y\Big|_{B_x}\simeq _{g_h}T^*Y|_{B_x}\simeq _{g_h}TY|_{B_x}\) over a small
  neighborhood \(B_x\) of \(x\) in \(Y_c\), where \(\simeq_{g_h}\)
  denote isomorphisms induced by the metric \(g_h\). Take a sequence
  \(\{y_i\}_i\subset B_x\) such that \(d\ff_h(y_i)\neq 0\) and \(y_i\to x\). (Such a
  sequence exists by Aronszajn's theorem.) Given \(\eta\neq 0\in
  (\bbR^3)^*\), use the same notation to denote the
  corresponding element in \(T^*Y_{y_i}\) or \(T^*Y_c\) under the
  aforementioned trivialization. Passing to a subsequence if
  necessary, we choose \(y_i\) to approach \(x\)  from the direction
  \(\beta \neq 0\in T_xY\simeq \bbR^3\). Let \(
  \rmh_i(y):=h_i(y_i)\delta _i(y, y_i)\), where \(h_i(y_i)\) is defined
  as in p.647 of \cite{H}, where \(\delta _i(y, y_i)\) are smooth compactly
  supported functions approximating the Dirac \(\delta \)-function
  \(\delta  (x,y)\) in the sense of distributions. Then a computation
  similar to that in Section 2.3 of \cite{H} shows that 
  \[
  \lim_{i\to\infty}  ev_{0, x}(g_h)(\rmh_i)=\lim_{i\to\infty}  \int_Y\big(d_1d_yG_{g_h}(y_i,y)\big)(\delta
    _{\rmh_i}*_{g_h})d_y\ff_{h}(y)=R_\beta (\eta),
    \]
  where \(R_\beta \) is the isomorphism defined in \cite{H}, and
  \(d_1\) denotes exterior derivative with respect to the first
  variable of the Green's function. 
\end{itemize}
\epf

\section{Some properties of the Seiberg-Witten solutions}

\subsection{Vortex solutions and the case when \(Y=\bbR\times\Sigma \)}\label{sec:vor}

\begin{prop}
Let  \(Y=\bbR_t\times \Sigma
 \) with the product metric, \(\ff=t\), 
 where \((\Sigma , \omega_\Sigma , j)\) is a
 compact K\"ahler surface 
 Let \(\grs_d\) be the \(\Spin^c\) structure
 on \(Y\) of degree \(d=d_\grs\), and let \(w=\pi _2^*\ul{w}\), where
 \(\ul{w}\) is a closed 2-form on \(\Sigma \) such that  \[
   \int_\Sigma
   \ul{w}=8\pi (d-\rmg+1), \]
 \(\rmg\) being the genus of \(\Sigma \). 
Then \(\forall r\geq 1\), there is a 1-1 map from \(\scrZ_{r, w}(\bbR_t\times \Sigma ,
 \grs_d, t)\) to \(\sym^d\Sigma \). 
 Here \(\sym^d\Sigma \) is
 defined to be \(\emptyset\) when \(d<0\), and \(\sym^0\Sigma \)
 consists of a point. 
\end{prop}

\pf Let \(\bbS\) denote the 
 spinor bundle corresponding to \(\grs_d\). Then \(\rho (dt)\) splits \(\bbS\) into a direct sum
 of eigen-subbundles \(E\), \(E\otimes K^{-1}\) corresponding to
 eigenvalues \(-i, i\) respectively:
 \[\bbS=E\oplus E\otimes K^{-1},\]
 where \(K^{-1}=\pi _2^*T\Sigma \), and \(E\simeq\pi _2^*E_\Sigma \), \(E_\Sigma \) being a complex line
 bundle over degree \(d\) over \(\Sigma \).
Write \[\Psi =2^{-1/2}r^{1/2} (\alpha , \beta )\in \Gamma  (Y; E\oplus E\otimes
  K^{-1}).\]
On the other hand, 
 Noting that a hermitian
 connection \(A^E\) on \(E\) induces a \(\Spin^c\) connection \(A\) on \(\bbS\) and
 vice versa, we will also use \((A^E, (\alpha , \beta ))\) to specify
 a configuration. 
 
Choose a reference connection \(A_0^E\) on \(E\) to be of the form \(A_0^E=\pi _2^*B_{\ul{w}}^E\), where
 \(B_w^E\) is a hermitian connection on \(E_\Sigma \) with
 \[F_{B_{\ul{w}}^E}=-i\ul{w}/4- F_{A^K}/2,\]  
 where \(A^K\) denotes the Levi-Civita curvature on the anti-canonical
         bundle \(K^{-1}=T^{0,1}\Sigma \). 
 Write
 \(a^E:=A^E-A^E_0=\rma_t(t, z)\, dt +a_z(t, z)\),
 where \(t\in \bbR\), \(z\in(\Sigma . j)\), \(\rma_t\) is an imaginery-valued
 function on \(Y\), and for each fixed \(t\), \(a_z(t, \cdot)\) is an
 imaginery-valued 1-form on \(\Sigma \),  Then a
 configuration \((A^E,(\alpha , \beta ))\) is \(l\)-admissible iff
 \(a_z\in \hat{L}^2_l(Y; \pi _2^*T^*\Sigma )\), \(\rma_t\in
 L^2_l(Y; i\bbR))\). Let \(\rma_z(t, \cdot)\in \Gamma (\Sigma ;
 \bbC)\) be defined by \(a_z(t, \cdot)=\rma_z(t, \cdot )dz+\bar{\rma}_z(t. \cdot)
 d\bar{z}\). Let \(B^E_z(t, \cdot)\) denote the connection on
 \(E\Big|_{\{t\}\times \Sigma }\simeq E_\Sigma \) given by
 \(B_w^E+a_z(t)\), and let \(B^{'E}_z(t, \cdot)\) denote the connection on
 \((E\otimes K^{-1})\Big|_{\{t\}\times \Sigma }\simeq E_\Sigma \otimes
 T\Sigma \) induced from  \(B^E_z(t, \cdot)\) and the Levi-Civita
 connection.
 
With such choices,  the Seiberg-Witten equation \(\grF_{\mu _{r,w}}(A, \Psi )=0\) takes the
 following form:
 \begin{eqnarray}\label{swr-1}
 &  \ds *_{g_\Sigma }d_z a_z+\frac{ir}{4}\big(1-|\alpha |^2+|\beta |^2)=0;\\
    &\ds \partial_t \rma_z-2\partial_z\rma_t=\frac{ir}{2}\beta \bar{\alpha
   }; \label{swr-2}\\
   & 2\bar{\partial}_{B^E_z}\alpha -(\partial_{t}+\rma_t) \beta =0;\label{swr-3}\\
   & 2 \partial_{B^{E'}_z}\beta +(\partial_{t}+\rma_t) \alpha =0,\label{swr-4}
 \end{eqnarray}
 where \(*_{g_\Sigma }\) denotes the two
 dimensional Hodge dual respect to the K\"ahler metric on \(\Sigma \),
  and \(d_z\) denotes the 2-dimensional exterior derivative in the
  \(z\)-variable.

 To proceed, note that any configuration \((A^E, (\alpha , \beta ))\)
 may be bring to one with \[\rma_z=0\]
by integrating along \(t\). (A configuration satisfying the above
equation is said to be in  a {\em temporal gauge}.)

\begin{lemma}
  Let \(l>1\) be an integer. 
  Then \(\forall r\geq 1\), any \(l\)-admissible
  solution  \((A^E, (\alpha , \beta ))\)  to  (\ref{swr-1}-\ref{swr-4}) in a temporal gauge
 satisfies  \(\beta \equiv
  0\). 
\end{lemma}
\pf
Combining the admissibility condition on \((A^E, (\alpha , \beta
))\), Sobolev embedding, and a Weitzenb\"ock formula, the
Seiberg-Witten equation  (\ref{swr-1}-\ref{swr-4}) implies
\[\begin{split}
  & \Big\langle \beta ,
  \big(\nabla_{A^{'E}}^*\nabla_{A^{'E}}+\frac{r}{4}(1+|\alpha
  |^2+|\beta |^2)\big)\beta \Big\rangle_{L^2}\\
&\quad =  \|\nabla_{A^{'E}}\beta \|^2_{L^2}+\frac{r}{4}\int_Y\big(1+|\alpha
  |^2+|\beta |^2\big)|\beta |^2=0, 
\end{split}
\]
This implies that \(\beta \equiv 0\) if \(r>0\).
\epf

Now, 
set \(\beta =0\), \(\rma_t=0\) in
(\ref{swr-1}-\ref{swr-4}). 
This implies that \(\partial_t a_z=0\), \(\partial_t\alpha =0\), that is,
\((a_z, \alpha )=\pi _2^*(\ul{a}, \ul{\alpha })\), where \(\ul{a}\) is
a connection on \(E_\Sigma \), and \(\ul{\alpha }\) is a section of
\(E_\Sigma \). Moreover, \((\ul{a}, \ul{\alpha })\) satisfies
\begin{equation}\label{eq:vor}
\grV_{r, d}(\ul{a}, \ul{\alpha}):=\left( \begin{array}{c}
             *_{g_\Sigma }d_z\ul{a}+\frac{i}{2}r(1-|\ul{\alpha }|^2)\\
    \bar{\partial}_{\ul{a}}\ul{\alpha }
           \end{array}
         \right)=0. 
\end{equation}
Equivalently, \((B_d^E+\ul{a}, r^{1/2}\ul{\alpha })\)
satisfies the {\em vortex equation} on \(\Sigma \), as defined  in \cite{G}
Equation (2), with the parameter \(\tau=r+\rmc d_\grs\). Here,  \(\rmc:=\ds \frac{8\pi
 }{\int_\Sigma \omega_\Sigma }\), and \(B_d^E\) is a connection on
 \(E_\Sigma \) with \(\ds F_{B_d^E}=-i\frac{\rmc }{4}d_\grs\, \omega_\Sigma
 \). 

The vortex equation is invariant under the action of
\(C^\infty(\Sigma ; U(1))\): given \(u\in \scrG_\Sigma :=C^\infty(\Sigma ; U(1))\),
\[
  u\cdot (\ul{a}, \ul{\alpha }):=(\ul{a}-u^{-1}du, u\cdot\ul{\alpha }).
\]
We denote by \(\scrV_{r,d}(\Sigma )\) the {\em moduli space
  of vortex solutions},
\begin{equation}\label{def:mod-vor}
  \scrV_{r, d}(\Sigma ):=\grV_{r,d}^{-1}(0)/\scrG_\Sigma .
\end{equation}

Given a pair \((\ul{a}, \ul{\alpha })\in i\Omega ^1(\Sigma )\times
\Gamma (E_\Sigma )\), we call
the Seiberg-Witten configuration \((A^E, (\alpha , \beta ))=(\pi _2^*(B_{\ul{w}}^E+\ul{a}), (\pi
_2^*\ul{\alpha }, 0))=:\jmath (\ul{a}, \ul{\alpha })\) the {\em pullback} of \((\ul{a}, \ul{\alpha
})\). We saw that \(\jmath\) defines a 1-1 map from the space of
solutions to the vortex equation (\ref{eq:vor}) to the space of
Seiberg-Witten solutions in temporal gauge. 

Meanwhile, observe tha two Seiberg-Witten configurations in temporal
gauge (in particular, pullback configurations) are
gauge-equivalent iff they are related by a gauge action by \(\pi
_2^*u\) for certain \(u\in C^\infty(\Sigma ; U(1))\), and
\[
  (\pi _2^*u)\cdot  \jmath (\ul{a}, \ul{\alpha })=\jmath \big( u\cdot (\ul{a}, \ul{\alpha })\big).
  \]
Thus, \(\jmath\) defines a 1-1 map from \(\scrV_{r, d}(\Sigma )\) to  \(\scrZ_{r, w}(\bbR_t\times \Sigma ,
 \grs_d, t)\)

By Theorem 1 of \cite{G}, \(\scrV_{r,d}(\Sigma )\) is diffeomorphic to
\(\op{Sym}^d\Sigma \) when \(\tau >\rmc d\), namely, when
\(r>0\). Moreover, it is endowed with a symplectic structure induced
from its embedding as a symplectic quotient in \(i\Omega ^1(\Sigma )\times
\Gamma (E_\Sigma )\), the latter being equipped with a natural
symplectic form (cf. \cite{G} p.91).
\epf

\subsection{Some properties of vortex solutions}

We shall need the following well-known property of vortex solutions.

\subsubsection{Pointwise estimates}
\begin{lemma}\label{lem:vor-Phi}
  Given \(  (\ul{a}, \ul{\alpha })\in \scrV_{r, d}(\Sigma )\),
  \[
    \|\ul{\alpha }\|^2_{L^\infty}\leq  1+r^{-1}\|\rms_-\|_{L^\infty}, 
  \]
  where \(\rms\) is the scalar curvature, \(\rms_-(z):=\op{max}
  (-\rms(z), 0)\).
\end{lemma}

A proof can be given along the line of the proof of Lemma
\ref{lem:Linf} below.

\subsubsection{Local structure of \(\scrV_{r, d}(\Sigma )\)}

\subsubsection{The symplectic structure on \(\scrV_{r, d}(\Sigma )\)}

\subsection{A priori estimates}

\subsubsection{An \(L^\infty \) bound on \(\Psi \) and \(F_A\).}   
Let \(\psi :=2^{1/2}r^{-1/2}\Psi \). We have the following standard
 $L^{\infty}$ bound on $\psi$ when $(A, \Psi )$ is an admissible
 solution to (\ref{eq:SW3}).

 Let \((Y, \grs)\) be a \(\Spin^c\) CCE, and let \(\rms\) denote the
scalar curvature. 
The constraint on the metric on \(Y\) implies that  $\|\max
(-\rms, 0)\|_{L^\infty(Y)}$ is well defined. Let \(\ff\) be a harmonic
function on \(Y\) satisfying Condition (2) of Definition
\ref{def:adm-f}. According to Proposition \ref{prop:adm-f},
\(\ds \|df\|_{L^\infty}\) is also finite. Fix an integer \(l>3\), and let \(w\in
\hat{\scrW}_{l,\grs}\). The admissibility condition on \(w\), together
with a version of Sobolev embedding, 
shows that \(\|w\|_{C^1}\) is
also finite. 

\begin{lemma}\label{lem:Linf}
Let \((Y, \grs)\), \(\ff\), \(w\) be as the above. Then any solution \((A, \Psi )\) to the Seiberg-Witten
equation \(\grF_{\mu _{r,w}}(A, \Psi )=0\) satisfies:
\begin{equation}\label{eq:Linf}
\|\psi\|^2_{L^\infty} \leq \|d\ff\|_{L^\infty}+z'r^{-1},
\end{equation}
where \(z'\) is  a positive constant depending only
on the \(L^\infty\) bounds on \(\rms\), \(w\) mentioned previously.

Via the curvature equation in (\ref{eq:SW3}), this gives
an $L^{\infty}$-bound for $F_{A^t}$:
\begin{equation}\label{F-Linf}
  \|F_{A^t}\|^2_{L^\infty} \leq 2r\|d\ff\|_{L^\infty}+z'',
\end{equation}
where \(z''\) is  a positive constant depending only
on the \(L^\infty\) bounds on \(\rms\), \(w\). 
\end{lemma}

\pf The Dirac equation in (\ref{eq:SW3}) together with a Weitzenb\"ock
formula gives: 
\begin{equation}
\slp_A \slp_A \psi=\nabla_A^*\nabla_A\psi+\frac{\rms}{4} \psi+\frac{\rho(F_A)}{2} \psi =0.
\end{equation}
Taking pointwise inner product of the preceding equation with $\psi$,
and using the curvature equation in (\ref{eq:SW3}),
we have
\begin{equation}\label{ineq:psi}\begin{split}
\frac{1}{2}d^*d|\psi|^2 & +|\nabla_A\psi|^2+\frac{r}{4}
|\psi|^2(|\psi|^2-r^{-1}|\mu _{r,w}|)+\frac{\rms}{4}|\psi|^2=0
\end{split}
\end{equation}
The smooth function \(|\psi|^2\) must have a maximum at a certain
point \(x_M\in Y\), or it is
bounded by one of \(2r^{-1}\|\Phi_\pm\|_{L^\infty(\Sigma _\pm, \bbS_{\Sigma _\pm})}\), where \((B_\pm,
\Phi_\pm)=\Pi _{\pm\infty}(A, \Psi )\). In the former
case, consider the previous inequality at \(x_M\) and rearranging to
get
\[
|\psi(x_M)|^2\big(|\psi( x_M)|^2-|df(x_M)|^2-z'r^{-1}\big)\leq 0
\]
where \(z'\) is  a positive constant depending only
on the \(L^\infty\) bounds on \(\rms\), \(w\) mentioned
previously. This leads directly to (\ref{eq:Linf}). 
In the latter case, invoke
Lemma \ref{lem:vor-Phi} and the fact that \(\ds \|d\ff\|_{L^\infty}\geq
1\).
\epf

\subsubsection{A pointwise bound for \(|\beta |^2\).}

Let \((Y, \grs)\) be a \(\Spin^c\) CCE, and \(\bbS\) be the
corresponding \(\Spin^c\) bundle. 
Let \(\ff\) be an
admissible function on \(Y\).
Let \(Z_{\ff} \subset Y_c\) consists of the critical points of
\(\ff\). Then over \(Y':=Y\backslash Z_{\ff}\), let \(K^{-1}\) be the
subbundle \(\ker (d\ff)\subset TY|_{Y'}\), endowed with the complex
structure given by the Clifford action of \(\rho
(d\ff)/|d\ff|\). Let \(A^K\) be the connection on \(K^{-1}\) induced
from the Levi-Civita connection. 
Split
\begin{equation}\label{def:S-split}
  \bbS|_{Y'}=E\oplus  E\otimes K^{-1}
\end{equation}
as a direct sum of eigenbundle of \(\rho (d\ff)\), where \(E\) is the
eigenbundle with eigenvalue \(-i|d\ff|\). Given a \(\Spin^c\)
connection \(A\) on \(\bbS\), denote by \(A^E\), \(A^{'E}\)
respectively the induced connection on \(E\), \(E\otimes
K^{-1}\). For simplicity, we shall use \(\nabla_A\) to denote
covariant derivatives with respect to any connection induced from
\(A\) and the Levi-Civita connection. For example, \(\nabla_A\alpha
=\nabla_{A^E}\alpha \); \(\nabla_A\beta =\nabla_{A^{'E}}\beta \). 
Given \(\Psi \in \Gamma (\bbS)\), write
\[
  \Psi |_{Y'}=2^{-1/2}r^{1/2}(\alpha , \beta )
\]
according to the splitting (\ref{def:S-split}).

Let \(\ts\) be the function on \(Y'\) defined as follows: 
Suppose that \(|\nu|^{-1}(0)\neq\emptyset\). Let \(\sigma(\cdot)\) denote the distance function to
\((d\ff)^{-1}(0)\), and set 
\[
\tilde{\sigma}:=(1-\chi(\sigma))\, \sigma+\chi(\sigma).
\]
When \(\ff\) has no critical points, let \(\sigma =\infty\) and
\(\td{\sigma }=1\). 

Let \(Y'_\delta:=\{x|\, \sigma (x)\geq \delta \} \subset Y\).

\begin{lemma}\label{lem:b0}
  Let \((Y, \grs)\), \(\ff\) be as the above, and let \(w\in
  \hat{\scrW}\). Let \((A, \Psi )\) be an admissible solution to
  \(\grF_{\mu _{r, w}}(A, \Psi )=0\).

There exist positive constants \(\sO\geq 8\), \(c\), \(c'\)
\( \zeta_0, \zeta '_0\geq 1\)  that depend only on
the metric, \(\ff\), and \(w\), 
such the following hold: Suppose \(r>1\), \(\delta >0\) are such that
\(\delta\geq\textsc{o}r^{-1/3}\), then 
\begin{equation}\label{ineq:beta0}
\begin{split}
|\beta|^2& \leq 2 c \, \ts^{-3}r^{-1} (
|d\ff|-|\alpha|^2)+\zeta _0\, \ts^{-5} r^{-2};\\
|\beta|^2& \leq 2 c' \ts^{-3}r^{-1} (
|d\ff|-|\psi|^2)+\zeta '_0\, \ts^{-5} r^{-2}
\end{split}
\end{equation}
on  \(Y_\delta\). 
\end{lemma}
\pf
This follows from a straightforward adaption of the proof of Proposition 5.5 of
\cite{L3}. 
\epf

\begin{lemma}\label{lem:est-1st-der}
There exist positive constants \(r_1\), \(\zeta _O\), \(\zeta', \zeta''\), that are independent of
\(r\) and \((A, \Psi )\), with the following significance: 
Let \(\delta _0'=\zeta _O r^{-1/3}\). For any \(r>r_1\), one has: 
\[\begin{split}
|\nabla_A\underline{\alpha}|^2+r\ts^2|\nabla_A\beta|^2&\leq \zeta'r
\varpi+\zeta''\ts^{-2} \quad \text{over \(Y_{\delta'_0}\)}.
\end{split}\]
\end{lemma}
\pf
This follows from  straightforward adaption of the proof of
Proposition 5.9 in \cite{L3}. The argument is much simpler here, since
instead of the complicated curvature estimates in the 4-dimensional
setting of \cite{L3}, in the
3-dimensional case, the required curvature estimate follows readily
from Lemma \ref{lem:Linf}.
\epf

\subsection{Asymptotic behaviors of admissible Seiberg-Witten
   solutions}

 \subsubsection{End-point maps from Seiberg-Witten moduli space to vortex
   moduli spaces.} Observe that if \((A, \Psi )\) is an (\(l\)-) admissible
 solution to \(\grF_{\mu _r}(A, \Psi )=0\), Then \((B^E_\pm, \Phi _\pm):=\Pi _{\pm\infty}(A,
 \Psi )\in \Conn (E_{\Sigma _\pm})\times\Gamma (E_{\Sigma _\pm})\)
 must be a vortex solution. More precisely, write \(B^E_\pm= B^E_{\pm,
   0}+\ul{a}_\pm\), \(\Phi =r^{1/2}2^{-1/2}\ul{\alpha }_\pm\), then \((\ul{a}_\pm,
 \ul{\alpha }_\pm)\) must satisfy \(\grV_{r, d_\pm}(\ul{a}_\pm,
 \ul{\alpha }_\pm)=0\). This induces end-point maps
 \[
   \Pi _{\pm\infty}\co \scrZ_{r,w}(Y, \grs, \ff)\to \scrV_{r,
     d_\pm}(\Sigma _\pm). 
   \]

\subsubsection{Exponential Decay of \(|\beta |^2\) and \(|\nabla_A\beta |^2\).}

\begin{lemma}\label{lem:b1}
  Let \((Y, \grs)\), \(\ff\) be as the above, and let \(w\in
  \hat{\scrW}\). Let \((A, \Psi )\) be an admissible solution to
  \(\grF_{\mu _{r, w}}(A, \Psi )=0\). Let \(\epsilon >0\) be as in
  (\ref{def:e}). Then there exist constants \(r_0\geq 1\), \(C>0\) depending only on the metric, \(w\), and
  \(d\ff\), such that \(\forall r\geq r_0\), the following holds:  Over \(\scrE_\pm[Y]\), we
  have the pointwise bound
  \[
    |\beta |^2 \leq C r^{-1}e^{-2\epsilon \tilde{t}}, 
  \]
  where \(C>0\) is a constant depending only on the metric, \(\epsilon
  \), \(w\), and
  \(d\ff\).  
\end{lemma}
\pf
Take pointwise inner product of the equation \(\slp_A^2\psi =0\) with
\(\beta \) and \(\alpha \) respectively to get the analogs of
Equations (2.3) and (2.4) of \cite{T}:
\begin{eqnarray}
\big(\frac{d^*d}{2}+r\frac{|\psi |^2+|d\ff|}{4}\big) |\beta|^2
                   +|\nabla_A\beta|^2
                   \leq \big(\zeta_1|b|\, |\nabla_A\alpha|+ \zeta
                   _1'|\nabla b|\, |\alpha | \big)|\beta|, 
\label{ineq:beta}\\
\frac{d^*d}{2} |\alpha|^2
                      +|\nabla_A\alpha|^2-\frac{r}{4}(|d\ff|-|\psi
                      |^2) |\alpha|^2 \leq
                      \big(\zeta_2|b|\, |\nabla_A\beta |+ \zeta
                   _2'|\nabla b|\, |\beta | \big)|\alpha |,
\label{ineq:alpha}
\end{eqnarray}
where \(b\) arises from
\(\nabla (d\ff)\),  and by Proposition \ref{prop:adm-f}, we have 
\begin{equation}\label{bdd:b}
  |b|+|\nabla b|\leq\zeta'_0\, e^{-\epsilon \tilde{t}}
\end{equation}
on
\(\scrE_\pm[Y]\). In the above as well as for the rest of this proof, the positive constants
\(\zeta _i, \zeta _i'\) depend only on
the metric, \(d\ff\), and \(w\).

Using Proposition \ref{prop:adm-f}
again, we may choose \(R>0\) such that \(\frac{1}{2}\leq |d\ff|\leq 2
\) over \(\scrE^R_\pm[Y]\). Assume also that \(r>1\) is much larger than
the \(L^\infty\) bound of \(w\) and \(\rms\). 
Then applying a triangular inequality 
to (\ref{ineq:beta}) and rearranging, one has
\begin{equation}\label{b1}
  \begin{split}
    & \big(\frac{d^*d}{2}+r\frac{|\psi |^2+|d\ff|}{8}\big) |\beta|^2
                   +|\nabla_A\beta|^2
   \leq  \zeta _3 r^{-1}e^{-2\epsilon \tilde{t}}|\nabla_A\alpha|^2
\end{split}
\end{equation}
over \(\scrE^R_\pm[Y]\). 
Meanwhile, write \(\varpi:=|d\ff|-|\alpha |^2\), and note that by
Proposition \ref{prop:adm-f},
\[
 \Big| d^*d|d\ff|\Big|\leq z_0 \, e^{-\epsilon \tilde{t}}\quad \text{over
   \(\scrE^R_\pm[Y]\),}
\]
where \(z_0>0\) is a constant depending only on the metric on
\(Y\). Combine the preceding inequality with (\ref{ineq:alpha}) as
well as Lemma \ref{lem:Linf} to get:
\begin{equation}\label{a1}
  \begin{split}
&\frac{d^*d}{2} (-\varpi) +|\nabla_A\alpha|^2+\frac{ r|\alpha |^2}{8}(-\varpi+|\beta |^2)\\
&\qquad \quad \leq e^{-2\epsilon \tilde{t}}\big(\zeta_3'|\nabla_A\beta|^2+\zeta_3'')  \quad
\text{over    \(\scrE^R_\pm[Y]\).}
\end{split}
\end{equation}
Adding \(\zeta _4r^{-1}e^{-2\epsilon \tilde{t}}\) times (\ref{a1}) to
(\ref{b1}) for an appropriately chosen constant \(\zeta _4>0\),  we
have for \(\rmu:=|\beta |^2-\zeta _4r^{-1}e^{-2\epsilon \tilde{t}}\varpi\):
\[
  \big(\frac{d^*d}{2}+r\frac{|d\ff|}{8}\big) \,\rmu\leq \zeta
  _5 r^{-1}e^{-2\epsilon \tilde{t}} \quad
\text{over    \(\scrE^R_\pm[Y]\).}
  \]
  Combine this with the fact that
  \begin{equation}\label{exp-comparison}
    d^*de^{-2\epsilon \tilde{t}} =-4\epsilon ^2e^{-2\epsilon
      \tilde{t}} \quad
\text{over    \(\scrE_\pm[Y]\)}
\end{equation}
as well as Lemmas \ref{lem:Linf}, \ref{lem:b0}, one may find a constant \(C'>0\) depending only on the metric, \(\epsilon
  \), \(w\), and \(d\ff\), such that
  \begin{equation}\label{bb}
    \begin{split}
   &  \big(\frac{d^*d}{2}+r\frac{|d\ff|}{8}\big)\,(\rmu-C'r^{-1}e^{-2\epsilon
      \tilde{t}}) <0 \quad
    \text{over    \(\scrE^R_\pm[Y]\)};\\
    & \big(\rmu-C'r^{-1}e^{-2\epsilon
      \tilde{t}}\big)\Big|_{\partial  \scrE^R_\pm[Y]}<0;\\
   & \Pi _{\pm\infty}\big(\rmu-C'r^{-1}e^{-2\epsilon
      \tilde{t}}\big)=0.
  \end{split}
\end{equation}
Suppose that there is an \(x\in \scrE^R_\pm[Y]\) where \(\rmv:=\rmu-C'r^{-1}e^{-2\epsilon
      \tilde{t}}>0\). Then \(\rmv\) attains a positive maximum in the
    interior of \(\scrE^R_\pm[Y]\). However, at such a maximum point,
    the left hand side in the first line of (\ref{bb}) is positive,
    which contradicts (\ref{bb}). Thus, \(\rmu\leq 0\) over
    \(\scrE^R_\pm[Y]\), which implies via Lemma \ref{lem:Linf} that
    \[
      |\beta |^2\leq C r^{-1}e^{-2\epsilon\tilde{t}}\quad \text{over    \(\scrE^R_\pm[Y]\).}
    \]
    Since \(\overline{\scrE_\pm[Y]\backslash\scrE^R_\pm[Y]}\) is
    compact, enlarging the value of the constant \(C\) if necessary,
    we arrive at the coclusion of the lemma.
    \epf

    \begin{lemma}\label{lem:b2}
  Let \((Y, \grs)\), \(\ff\) \(w\), \((A, \Psi )\) and  \(\epsilon >0\) be as in
  Lemma \ref{lem:b1}. Then there exist constants \(r_0\geq 1\), \(C'>0\) depending only on the metric, \(w\), and
  \(d\ff\), such that \(\forall r\geq r_0\), the following holds:  Over \(\scrE_\pm[Y]\), we
  have the pointwise bound
  \[
    |\nabla_A\beta |^2 \leq C r^{-1}e^{-2\epsilon \tilde{t}}, 
  \]
  where \(C>0\) is a constant depending only on the metric, \(\epsilon
  \), \(w\), and
  \(d\ff\).  
\end{lemma}
\pf Let \(R\) and \(r\) be sufficiently large positive numbers as in
the proof of the previous lemma. 
With (\ref{F-Linf}) in place, argue as in the proof of Proposition 2.8
in \cite{T} using this bound, (\ref{bdd:b}) and Lemma \ref{lem:Linf}
to  get:
\begin{equation}\label{DE:db}
  \begin{split}
 &  \big(\frac{d^*d}{2}+r\frac{|\psi |^2+|d\ff|}{4}\big) |\nabla_A\beta|^2
                   +|\nabla_A\nabla_A\beta|^2\\
&\quad                 \leq  \zeta _0 r |\nabla_A\beta|^2+
                   r^{-1}e^{-2\epsilon \tilde{t}}\big(\zeta_1
                   |\nabla_A\alpha|^2+ \zeta_2|\nabla_A\nabla _A\alpha
                   |^2+\zeta _3 \big) \quad \text{over    \(\scrE^R_\pm[Y]\)}
                 \end{split}
               \end{equation}
and
\begin{equation}\label{DE:da}
  \begin{split}
 &  \frac{d^*d}{2}|\nabla_A\alpha |^2
                   +|\nabla_A\nabla_A\alpha |^2\\
&\quad                  \leq  \zeta '_0 r |\nabla_A\alpha |^2+
                   r^{-1}e^{-2\epsilon \tilde{t}}\big(\zeta'_1
                   |\nabla_A\beta |^2+ \zeta'_2|\nabla_A\nabla _A\beta |^2
                   |+\zeta '_3 |\beta |^2\big) \quad \text{over    \(\scrE^R_\pm[Y]\).}
                 \end{split}
               \end{equation}

Adding \(C_1'r^{-1}e^{-2\epsilon \tilde{t}}\) times (\ref{DE:da}) to
(\ref{DE:db}) for an appropriately chosen constant \(C_1'\), we have
for \(\rmu_1:=|\nabla_A\beta|^2-C_1'r^{-1}e^{-2\epsilon
  \tilde{t}}|\nabla_A\alpha |^2\):
  \[
  \frac{d^*d}{2}\rmu_1\leq  \zeta '_4 r |\nabla_A\beta |^2+\zeta _5'r^{-1}e^{-2\epsilon \tilde{t}} |\nabla_A\alpha |^2
                    \quad \text{over    \(\scrE^R_\pm[Y]\). }\]
Adding \(\zeta _4'r\) times (\ref{b1}) to the preceding inequality, we
have:
\[
  \frac{d^*d}{2}(\rmu_1+\zeta _4'r|\beta |^2) \leq  \zeta _6'r^{-1}e^{-2\epsilon \tilde{t}} |\nabla_A\alpha |^2
                    \quad \text{over    \(\scrE^R_\pm[Y]\). }
                  \]
Using (\ref{exp-comparison}) and Lemma \ref{lem:est-1st-der}, we may
find another positive constant \(\zeta _7\), such that with
\(\rmv':=\rmu_1+\zeta _4'r|\beta |^2-\zeta _7e^{-2\epsilon
  \tilde{t}}\),
\begin{equation}
    \begin{split}
   &  \frac{d^*d}{2}\,(\rmv') <0 \quad
    \text{over    \(\scrE^R_\pm[Y]\)};\\
    & \rmv'\Big|_{\partial  \scrE^R_\pm[Y]}<0;\\
   & \Pi _{\pm\infty}\rmv'=0.
  \end{split}
\end{equation}
A maximum principle type argumet as that in the proof of the previous
lemma then yields:
\[
  |\nabla_A\beta|^2-C_1'r^{-1}e^{-2\epsilon
  \tilde{t}}|\nabla_A\alpha |^2+\zeta _4'r|\beta |^2\leq \zeta _7e^{-2\epsilon
  \tilde{t}}. 
\]
A combination of the preceding inequality with Lemma
\ref{lem:est-1st-der} then leads to the conclusion of the lemma.
\epf

\subsubsection{An alternative parametrization of \(\scrE_\pm[Y]\) and
  reference pullback configurations}

We aim to show that an admissible Seiberg-Witten solution \((A;\Psi
)\) ``approaches the end-point vortex solutions \(\Pi _{\pm\infty}(A,
\Psi )\) exponentially''. To state this precisely, we shall
construct a reference configuration on \(Y\) from \(\Pi _{\pm\infty}(A,
\Psi )\), which approximate the pullback configurations on cylinders
defined in Section \ref{sec:vor} on the ends of \(Y\), then show that
the difference between \((A, \Psi )\) and this reference configuration
decays exponentially over \(\scrE_\pm[Y]\). This is done similarly to
what appears in Section 3.3 of \cite{L1}.


\end{document}